\documentclass{amsart}
\usepackage{amsmath, amscd, amssymb, fullpage}

\newtheorem{theorem}{Theorem}[section]

\newtheorem{corollary}[theorem]{Corollary}
\newtheorem{fact}[theorem]{Fact}
\newtheorem{conjecture}[theorem]{Conjecture}

\theoremstyle{definition}

\newtheorem{remark}[theorem]{Remark}
\newtheorem{definition}[theorem]{Definition}

\DeclareMathOperator{\Gal}{Gal(K^{\sep}/K)}
\DeclareMathOperator{\Stab}{Stab}
\DeclareMathOperator{\hhat}{\widehat{h}}
\newcommand{\tensor}{\otimes}

\def\Frac{\operatorname{Frac}}
\def\trdeg{\operatorname{trdeg}}

\def\alg{\operatorname{alg}}

\def\sep{\operatorname{sep}}

\def\dim{\operatorname{dim}}

\def\tor{\operatorname{tor}}

\begin{document}

\title{Towards the full Mordell-Lang conjecture for Drinfeld modules}
\author{Dragos Ghioca}
\address{Dragos Ghioca, Department of Mathematics \& Statistics, Hamilton Hall, Room 218, McMaster University, 1280 Main Street West, Hamilton, Ontario L8S 4K1, Canada}
\email{dghioca@math.mcmaster.ca}

\begin{abstract}
Let $\phi$ be a Drinfeld module of generic characteristic, and let $X$ be a sufficiently generic affine subvariety of $\mathbb{G}_a^g$. We show that the intersection of $X$ with a finite rank $\phi$-submodule of $\mathbb{G}_a^g$ is finite.
\end{abstract}

\maketitle

\section{Introduction}
\label{se:intro}

In \cite{McQuillan}, McQuillan proved the Mordell-Lang conjecture in its most general form.
\begin{theorem}[The full Mordell-Lang theorem]
\label{full-ml}
Let $G$ be a semi-abelian variety defined over a number field $K$. Let $X\subset G$ be a $K^{\alg}$-subvariety, and let $\Gamma\subset G(K^{\alg})$ be a finite rank group (i.e. $\Gamma$ lies in the divisible hull of a finitely generated subgroup of $G(K^{\alg})$). Then there exist algebraic subgroups $B_1,\dots,B_l$ of $G$ and there exist $\gamma_1,\dots,\gamma_l\in\Gamma$ such that 
$$X(K^{\alg})\cap\Gamma = \bigcup_{i=1}^l \left(\gamma_i + B_i(K^{\alg})\right)\cap\Gamma.$$
\end{theorem}

We note that in Theorem~\ref{full-ml} if $X$ does not contain any translate of a positive dimensional algebraic subgroup of $G$, then the full Mordell-Lang theorem says that $X(K^{\alg})\cap\Gamma$ is finite.
Also, a particular case of the full Mordell-Lang theorem (in the case $\Gamma$ is the torsion subgroup $G_{\tor}$ of $G$) is the Manin-Mumford theorem, which was first proved by Raynaud \cite{Raynaud}. 

Faltings \cite{Faltings} proved the Mordell-Lang conjecture for finitely generated subgroups $\Gamma$ of abelian varieties $G$. His proof was extended by Vojta \cite{Vojta} to finitely generated subgroups of semi-abelian varieties $G$. Finally, McQuillan \cite{McQuillan} extended Vojta's result to finite rank subgroups $\Gamma$ of semi-abelian varieties $G$. Later, R\"{o}ssler \cite{Roessler} provided a simplified proof of McQuillan's extension in which he used uniformities for the intersection of translates of a fixed subvariety $X\subset G$ with the torsion subgroup of the semi-abelian variety $G$. Essentially, R\"{o}ssler showed that the full Mordell-Lang conjecture follows from the Mordell-Lang statement for finitely generated subgroups, combined with a uniform Manin-Mumford statement as proved by Hrushovski \cite{Hrushovski}. 

It is important to note that the exact translation of the Mordell-Lang conjecture to semi-abelian varieties in characteristic $p$ is false due to the presence of isotrivial varieties. However, Hrushovski \cite{Hru} saved the Mordell-Lang theorem for finitely generated subgroups of semi-abelian varieties in characteristic $p$ by treating isotrivial varieties as \emph{special}. The isotrivial case was treated by Rahim Moosa and the author in \cite{full-ml} where it was obtained a full Mordell-Lang statement for isorivial semi-abelian varieties in characteristic $p$. On the other hand, if we replace $G$ by a power $\mathbb{G}_a^g$ of the additive group scheme, then the exact translation of the Mordell-Lang conjecture either fails (in characteristic $0$) or it is trivially true (in characteristic $p$).

Inspired by the analogy between abelian varieties in characteristic $0$ and Drinfeld modules of generic characteristic, Denis \cite{Denis-ML} proposed that analogs of the Manin-Mumford and Mordell-Lang theorems hold for such Drinfeld modules $\phi$ acting on $\mathbb{G}_a^g$ (in characteristic $p$). Denis conjectures ask for describing the intersection of an affine subvariety $X\subset \mathbb{G}_a^g$ with a finite rank $\phi$-submodule $\Gamma$ of $\mathbb{G}_a^g$. Using methods of model theory, combined with some clever number theoretical arguments, Scanlon \cite{Scanlon} proved the Denis-Manin-Mumford conjecture. In \cite{IMRN}, the author proved the Denis-Mordell-Lang conjecture for finitely generated $\phi$-modules $\Gamma$ under two mild technical assumptions. In this paper, we extend our result from \cite{IMRN} to finite rank $\phi$-submodules $\Gamma$.

We also note that recently there have been significant progress on establishing additional links between classical diophantine results over number fields and similar statements for Drinfeld modules. The author proved in \cite{Math.Ann} an equidistribution statement for torsion points of a Drinfeld module, which is similar to the equidistribution statement established by Szpiro-Ullmo-Zhang \cite{Szpiro-Ullmo-Zhang} (which was later extended by Zhang \cite{Zhang} to a full proof of the famous Bogomolov conjecture). Also, Breuer \cite{Breuer} proved a special case of the Andr\'{e}-Oort conjecture for Drinfeld modules, while special cases of this conjecture in the classical case of a number field were proven by Edixhoven-Yafaev \cite{Edixhoven} and Yafaev \cite{Yafaev}. Bosser \cite{Bosser} proved a lower bound for linear forms in logarithms at an infinite place associated to a Drinfeld module (similar to the classical result obtained by Baker \cite{Baker} for usual logarithms, or by David \cite{David} for elliptic logarithms). Bosser's result was used by Thomas Tucker and the author in \cite{Tom} to establish certain equidistribution and integrality statements for Drinfeld modules. Moreover, Bosser's result is quite possibly true also for linear forms in logarithms at finite places for a Drinfeld module. Assuming this last statement, Thomas Tucker and the author proved in \cite{Tom-2} the analog of Siegel's theorem for finitely generated $\phi$-submodules. We believe that our present paper provides an additional proof of the fact that the Drinfeld modules represent the right arithmetic analog in characteristic $p$ for semi-abelian varieties in characteristic $0$.

The plan for our paper is as follows: in Section~\ref{se:statement} we provide the basic notation for our paper, while in Section~\ref{se:proof} we prove our main result (Theorem~\ref{T:full-ml-dr}).

\section{The Mordell-Lang theorem for Drinfeld modules}
\label{se:statement}

First we note that all subvarieties appearing in this paper are considered to be \emph{closed}. We define next the notion of a Drinfeld module.

Let $p$ be a prime and let $q$ be a power of $p$. 
Let $C$ be a projective non-singular curve defined over $\mathbb{F}_q$. Let $A$ be the ring of $\mathbb{F}_q$-valued functions defined on $C$, regular away from a fixed closed point $\infty\in C$. 
Let $K$ be a finitely generated field extension of the fraction field $\Frac(A)$ of $A$. We let $K^{\alg}$ be a fixed algebraic closure of $K$, and let $K^{\sep}$ be the separable closure of $K$ inside $K^{\alg}$.

We define the operator
$\tau$ as the Frobenius on $\mathbb{F}_q$, extended so that for every $x\in K^{\alg}$, we have $\tau(x)=x^q$. 
Then for every subfield $L\subset K^{\alg}$, we let $L\{\tau\}$ be the ring of polynomials in
$\tau$ with coefficients from $L$ (the addition is the usual addition, while the multiplication is given by the usual composition of functions).

Following Goss \cite{Goss}, we call a Drinfeld module of generic characteristic defined over $K$ a morphism $\phi:A\rightarrow
K\{\tau\}$ for which the 
coefficient of $\tau^0$ in $\phi_a$ is $a$ for
every $a\in A$, and there 
exists $a\in A$ 
such that $\phi_a\ne a\tau^0$. 
For the remainder of this paper, unless otherwise stated, $\phi:A\rightarrow K\{\tau\}$ is a Drinfeld module of generic characteristic.

A Drinfeld module $\psi:A\rightarrow K^{\alg}\{\tau\}$ is isomorphic to $\phi$ (over $K^{\alg}$) if there exists a nonzero $\gamma\in K^{\alg}$ such that for every $a\in A$, we have $\psi_a = \gamma^{-1}\phi_a\gamma$.

For every field extension $K\subset L$, the Drinfeld module $\phi$ induces an action on $\mathbb{G}_a(L)$ by $a*x:=\phi_a(x)$, for each $a\in A$. Let $g$ be a fixed positive integer. We extend diagonally the action of $\phi$ on $\mathbb{G}_a^g$. 

The subgroups of $\mathbb{G}_a^g(K^{\alg})$ invariant under the action of $\phi$ are called \emph{$\phi$-submodules}. For a $\phi$-submodule $\Gamma$, its full divisible hull is
$$\Gamma\tensor _A\Frac(A) := \left\{x\in \mathbb{G}_a^g(K)\mid \text{ there exists } 0\ne a\in A\text{ such that }\phi_a(x)\in\Gamma\right\}.$$
We define the \emph{rank} of a $\phi$-submodule $\Gamma\subset \mathbb{G}_a^g(K^{\alg})$ as $\dim_{\Frac(A)}\Gamma\tensor_A\Frac(A)$. 

\begin{definition}
An algebraic $\phi$-submodule of $\mathbb{G}_a^g$ is an algebraic subgroup of $\mathbb{G}_a^g$ invariant under $\phi$.
\end{definition}

Denis proposed in Conjecture $2$ of \cite{Denis-ML} the following problem.
\begin{conjecture}[The full Denis-Mordell-Lang conjecture]
\label{con:full-ml-dr}
Let $X\subset\mathbb{G}_a^g$ be an affine variety defined over $K^{\alg}$. Let $\Gamma$ be a finite rank $\phi$-submodule of $\mathbb{G}_a^g(K^{\alg})$. Then there exist algebraic $\phi$-submodules $B_1,\dots,B_l$ of $\mathbb{G}_a^g$ and there exist $\gamma_1,\dots,\gamma_l\in\Gamma$ such that 
$$X(K^{\alg})\cap\Gamma = \bigcup_{i=1}^l (\gamma_i+B_i(K^{\alg}))\cap\Gamma.$$
\end{conjecture}

Before stating our result, we need to introduce the following notion.

\begin{definition}
\label{D:modular transcendence degree}
We call the modular transcendence degree of $\phi$ the smallest integer $d\ge 1$ such that a Drinfeld module isomorphic to $\phi$ is defined over a field of transcendence degree $d$ over $\mathbb{F}_q$.
\end{definition}

In \cite{IMRN} (see Theorem $4.11$), the author proved the following result towards Conjecture~\ref{con:full-ml-dr}.

\begin{theorem}
\label{T:ml-dr}
With the above notation, assume in addition that the modular transcendence degree of $\phi$ is at least $2$. Let $X\subset\mathbb{G}_a^g$ be an affine subvariety defined over $K^{\alg}$ such that there is no positive dimensional algebraic subgroup of $\mathbb{G}_a^g$ whose translate lies inside $X$. Let $\Gamma$ be a finitely generated $\phi$-submodule of $\mathbb{G}_a^g(K^{\alg})$. Then $X(K^{\alg})\cap\Gamma$ is finite.
\end{theorem}

In Theorem~\ref{T:full-ml-dr} we extend the previous result to all finite rank $\phi$-submodules $\Gamma$.

\begin{remark}
\label{R:conditions}
We have two technical conditions in Theorem~\ref{T:ml-dr} which we will keep also in our extension from Theorem~\ref{T:full-ml-dr}. The condition that $\phi$ has modular transcendence degree at least equal to $2$ is a mild technical condition, however necessary due to the methods employed in \cite{IMRN}. The condition that $X$ does not contain any translate of a positive dimensional algebraic subgroup of $\mathbb{G}_a^g$ is satisfied by all sufficiently generic affine subvarieties $X$. 
\end{remark}

\section{Proof of our main result}
\label{se:proof}

We continue with the notation from Section~\ref{se:statement}.
We define the torsion submodule of $\phi$ as
$$\phi_{\tor}=\{x\in \mathbb{G}_a^g(K^{\alg})\mid\text{ there exists }a\in A\setminus\{0\}\text{ such that }\phi_a(x)=0\}.$$

Next we state our main result.
\begin{theorem}
\label{T:full-ml-dr}
Let $K$ be a finitely generated field of characteristic $p$ and let $g$ be a positive integer. Let $\phi:A\rightarrow K\{\tau\}$ be a Drinfeld module of generic characteristic. Assume the modular transcendence degree of $\phi$ is at least $2$.
Let $X\subset\mathbb{G}_a^g$ be an affine subvariety defined over $K^{\alg}$ such that there is no positive dimensional algebraic subgroup of $\mathbb{G}_a^g$ whose translate lies inside $X$. Let $\Gamma$ be a finitely generated $\phi$-submodule of $\mathbb{G}_a^g(K)$, and let $\Gamma':=\Gamma\tensor_A \Frac(A)$. Then $X(K^{\alg})\cap\Gamma'$ is finite.
\end{theorem}

In our proof of Theorem~\ref{T:full-ml-dr} we need a uniform version of Scanlon's result from \cite{Scanlon}. He proved the Manin-Mumford theorem for Drinfeld modules in the following form (see his Theorem $1$).
\begin{theorem}
\label{Scanlon-MM}
Let $\phi:A\rightarrow K\{\tau\}$ be a Drinfeld module and let $X\subset \mathbb{G}_a^g$ be an affine variety defined over $K^{\alg}$. Then there exist algebraic $\phi$-submodules $B_1,\dots,B_l$ of $\mathbb{G}_a^g$ and there exist $\gamma_1,\dots,\gamma_l\in\phi_{\tor}^g$ such that
$$X(K^{\alg})\cap\phi_{\tor}^g = \bigcup_{i=1}^l \left(\gamma_i+B_i(K^{\alg})\right)\cap\phi_{\tor}^g.$$
\end{theorem}
In Remark $19$ from \cite{Scanlon}, Scanlon notes that his Manin-Mumford theorem for Drinfeld modules holds uniformly in algebraic families of varieties, i.e. if $X$ varies inside an algebraic family of varieties, then there exists a uniform bound on the degrees of the Zariski closures of $X(K^{\alg})\cap\phi_{\tor}^g$. In particular, we obtain the following uniform statement for translates of $X$.

\begin{corollary}
\label{C:needed uniformity}
With the notation for $\phi$ and $X$ as in Theorem~\ref{Scanlon-MM}, assume in addition that $X$ contains no translate of a positive dimensional algebraic subgroup of $\mathbb{G}_a^g$. Then there exists a positive integer $N$ such that for every $x\in\mathbb{G}_a^g(K^{\alg})$, the set $\left(x+X(K^{\alg})\right)\cap\phi_{\tor}^g$ has at most $N$ elements.
\end{corollary}

\begin{proof}
Because $X$ contains no translate of a positive dimensional algebraic subgroup of $\mathbb{G}_a^g$, then for every $x\in\mathbb{G}_a^g(K^{\alg})$, the algebraic $\phi$-modules $B_i$ appearing in the intersection $\left(x+X(K^{\alg})\right)\cap\phi_{\tor}^g$ are all $0$-dimensional. In particular, the set $\left(x+X(K^{\alg})\right)\cap\phi_{\tor}^g$ is finite. Thus, using the uniformity obtained by Scanlon for his Manin-Mumford theorem, we conclude that the cardinality of $\left(x+X(K^{\alg})\right)\cap\phi_{\tor}^g$ is uniformly bounded by some positive integer $N$, independent of $x$.
\end{proof}

We will also use the following fact in the proof of our Theorem~\ref{T:full-ml-dr}.
\begin{fact}
\label{F:torsion}
Let $\phi:A\rightarrow K\{\tau\}$ be a Drinfeld module. Then for every positive integer $d$, there exist finitely many torsion points $x$ of $\phi$ such that $[K(x):K]\le d$.
\end{fact}

\begin{proof}
If $x\in\phi_{\tor}$, then the canonical height $\hhat(x)$ of $x$ (as defined in \cite{Denis} and \cite{Wan}) equals $0$. Also, as shown in \cite{Denis}, the difference between the canonical height and the usual Weil height is uniformly bounded on $K^{\alg}$. Actually, Denis \cite{Denis} proves this last statement under the hypothesis that $\trdeg_{\mathbb{F}_q}K=1$. However, his proof easily generalizes to fields $K$ of arbitrarily finite transcendence degree. For this we need the construction of a \emph{coherent good set of valuations} on $K$ as done in \cite{mw2} (see also the similar construction of heights from \cite{Wan}). Essentially, a coherent good set $U_K$ of valuations on $K$ is a set of defectless valuations satisfying a product formula on $K$ (for more details, we refer the reader to Sections $2$ and $3$ of \cite{mw2}). Then Fact~\ref{F:torsion} follows by noting that there are finitely many points of bounded Weil height and bounded degree over the field $K$ (using Northcott's theorem applied to the global function field $K$).

Moreover, Corollary $4.22$ of \cite{mw2} provides an effective upper bound on the size of the torsion of $\phi$ over any finite extension $L$ of $K$ in terms of $\phi$ and the number of places of $L$ lying above places in $U_K$ of bad reduction for $\phi$. Because for each field $L$ such that $[L:K]\le d$, and for each place $v\in U_K$, there are at most $d$ places $w$ of $L$ lying above $v$, we conclude that there exists an upper bound for the size of torsion of $\phi$ over all field extensions of degree at most $d$ over $K$ in terms of $\phi$, $d$ and the number of places in $U_K$ of bad reduction for $\phi$.
\end{proof}

\begin{proof}[Proof of Theorem~\ref{T:full-ml-dr}.]
Because $\phi$ has generic characteristic, $\Gamma'\subset\mathbb{G}_a^g(K^{\sep})$. Hence $X(K^{\alg})\cap\Gamma'=X(K^{\sep})\cap\Gamma'$.

We note that our theorem is equivalent with showing that if $X$ is a positive dimensional irreducible subvariety of $\mathbb{G}_a^g$ satisfying the hypothesis of Theorem~\ref{T:full-ml-dr}, then $X(K^{\sep})\cap\Gamma'$ is not Zariski dense in $X$. Moreover, at the expense of moding out through the stabilizer $\Stab(X)$ of $X$ (which is a finite group according to the hypothesis of Theorem~\ref{T:full-ml-dr}), we may assume $\Stab(X)$ is trivial. Note that after moding out through the (finite) group $\Stab(X)$, the variety $X$ is still positive dimensional and irreducible. Finally, at the expense of replacing $K$ by a finite extension, we may assume that $X$ is defined over $K$, and also $\Gamma\subset\mathbb{G}_a^g(K)$.

Our proof follows the argument from \cite{full-ml} (which in turn was inspired by the argument from \cite{Roessler}). We assume by contradiction that $X(K^{\alg})\cap\Gamma'$ is Zariski dense in $X$.

We claim that for every $x\in\Gamma'$ and for every $\sigma\in\Gal$, we have $\sigma(x)-x\in\phi_{\tor}^g$. Indeed, because $x\in\Gamma'$, then there exists $0\ne a\in A$ such that $\phi_a(x)\in\Gamma$. Because $\Gamma\subset\mathbb{G}_a^g(K)$ and $\phi$ is defined over $K$, then 
$$\phi_a(\sigma(x))=\sigma(\phi_a(x))=\phi_a(x).$$
Hence $\phi_a(\sigma(x)-x)=0$, as desired. Therefore, if $x\in X(K^{\sep})$, then 
$\sigma(x)-x\in (-x+X)\cap\phi_{\tor}^g$ (we also used that $X$ is defined over $K$, and so, $\sigma(x)\in X$). Using Corollary~\ref{C:needed uniformity}, we obtain an upper bound on the number of conjugates of $x$, which gives us an upper bound, say $N$, for $[K(x):K]$. Implicitly, we also get $[K(\sigma(x)):K]\le N$. Because $$K(\sigma(x)-x)\subset K(x,\sigma(x)),$$ we conclude $[K(\sigma(x)-x):K]\le N^2$. Using that $\sigma(x)-x\in\phi_{\tor}^g$, and using Fact~\ref{F:torsion}, we conclude that the set
\begin{equation}
\label{finite set}
\{\sigma(x)-x\mid x\in\Gamma'\text{ and }\sigma\in\Gal\}\text{ is finite.}
\end{equation}

Assuming that $X(K^{\sep})\cap\Gamma'$ is Zariski dense in $X$ (and also using that $X$ is irreducible), then either $X(K)\cap\Gamma'$ is Zariski dense in $X$, or $(X(K^{\sep})\setminus X(K))\cap\Gamma'$ is Zariski dense in $X$. If the former statement holds, then $X$ has a Zariski dense intersection with a finitely generated group, as $\Gamma'\cap\mathbb{G}_a^g(K)$ is a finite rank subgroup of the tame module $\mathbb{G}_a^g(K)$ (see Theorem $1$ of \cite{Wan}). By Theorem~\ref{T:ml-dr}, $X(K)\cap\left(\Gamma'\cap\mathbb{G}_a^g(K)\right)$ is finite and hence it cannot be dense in $X$ (because $X$ is assumed to be positive dimensional). Therefore the latter case of the above dichotomy should occur. Using \eqref{finite set}, we conclude that there exists a \emph{nonzero} torsion point $y\in\phi_{\tor}^g$ such that the set 
$$\{x\in X(K^{\sep})\cap\Gamma'\mid \sigma(x)-x=y \text{ for some $\sigma\in\Gal$} \}$$
is Zariski dense in $X$. Therefore, $y\in\Stab(X)$, contradicting the fact that $X$ has trivial stabilizer. This concludes the proof of Theorem~\ref{T:full-ml-dr}.
\end{proof}


\begin{thebibliography}{9}

\bibitem{Baker}
A. Baker, \emph{Transcendental number theory}. Cambridge University Press, London-New York, 1975. x+147 pp.

\bibitem{Bosser}
V. Bosser, \emph{Minorations de formes lin\'{e}aires de logarithmes pour les modules de Drinfeld}. (French) [Lower bounds of linear forms in logarithms for Drinfeld modules] J. Number Theory {\bf 75} (1999), no. 2, 279--323.

\bibitem{Breuer}
F. Breuer, \emph{The Andr\'{e}-Oort conjecture for products of Drinfeld modular curves}. J. reine angew. Math. {\bf 579} (2005), 115--144.

\bibitem{David}
S. David, \emph{Minorations de formes lin\'{e}aires de logarithmes elliptiques}. (French) [Lower bounds for linear forms in elliptic logarithms] M\'{e}m. Soc. Math. France (N.S.) {\bf 62} (1995), iv+143 pp.

\bibitem{Denis}
L. Denis, \emph{Hauteurs canoniques et modules de
Drinfeld}. (French) [Canonical 
heights and Drinfeld modules] Math. Ann. {\bf 294} (1992),
no. 2, 213--223.

\bibitem{Denis-ML}
L. Denis, \emph{Diophantine geometry on Drinfeld modules}. The arithmetic 
of function fields (Columbus, OH, 1991), 285--302, Ohio State Univ. Math. Res. 
Inst. Publ., 2, de Gruyter, Berlin, 1992.

\bibitem{Edixhoven}
B. Edixhoven and A. Yafaev, \emph{Subvarieties of Shimura type}. Ann. of Math. (2) {\bf 157} (2003), no. 2, 621--645.

\bibitem{Faltings}
G. Faltings, \emph{The general case of S. Lang's conjecture}. Barsotti Symposium 
in Algebraic Geometry (Abano Terme, 1991), 175-182, Perspect. Math., 15, 
Academic Press, San Diego, CA, 1994.

\bibitem{IMRN}
D. Ghioca, \emph{The Mordell-Lang Theorem for Drinfeld modules}, Internat. Math. Res. Notices, {\bf 53} (2005), 3273--3307.

\bibitem{Math.Ann}
D. Ghioca, \emph{Equidistribution for torsion points of a Drinfeld module}. Math. Ann. {\bf 336} (2006), 841--865.

\bibitem{mw2}
D. Ghioca, \emph{The Lehmer inequality and the Mordell-Weil theorem for Drinfeld modules}. J. Number Theory {\bf 122} (2007), 37--68.

\bibitem{full-ml}
D. Ghioca and R. Moosa, \emph{Division points on subvarieties of isotrivial semiabelian varieties}. to appear in International Mathematics Research Notices.

\bibitem{Tom}
D. Ghioca and T. J. Tucker, \emph{Equidistribution and integrality for Drinfeld modules}. submitted for publication, available online at http://www.math.mcmaster.ca/\~{ }dghioca/papers/papers.html.

\bibitem{Tom-2}
D. Ghioca and T. J. Tucker, \emph{Siegel's Theorem for Drinfeld modules}. preprint, available online at http://www.math.mcmaster.ca/\~{ }dghioca/papers/papers.html.

\bibitem{Goss}
D. Goss, \emph{Basic structures of function field
arithmetic}, Ergebnisse 
der Mathematik und ihrer Grenzgebiete (3) [Results in
Mathematics and Related
  Areas (3)], 35. Springer-Verlag, Berlin, 1996.

\bibitem{Hru}
E. Hrushovski, \emph{The Mordell-Lang conjecture for function fields},
  J. Amer. Math. Soc. \textbf{9} (1996), no.~3,
  667--690. 

\bibitem{Hrushovski}
E. Hrushovski, \emph{The Manin-Mumford conjecture and the model theory of difference fields}. Annals of Pure and Applied Logic {\bf 112} (2001), 43--115.

\bibitem{McQuillan}
M.~McQuillan, \emph{Division points on semi-abelian varieties}. Invent.
  Math. \textbf{120} (1995), no.~1, 143--159.

\bibitem{Raynaud}
M. Raynaud, \emph{Sous-vari\'{e}t\'{e}s d'une vari\'{e}t\'{e} ab\'{e}lienne et points de torsion}. (French) [Subvarieties of an abelian variety and torsion points]. Arithmetic and Geometry, Vol. I, 327-352, Prog. Math., {\bf 35}, Birkhauser, Boston, MA, 1983.

\bibitem{Roessler}
D. R\"{o}ssler, \emph{An afterthought on the generalized Mordell-Lang conjecture}. to appear in the Proceedings of the workshop on Model Theory, Algebraic and Analytic Geometry (Newton Institute, Cambridge, 11-15 July 2005).

\bibitem{Scanlon}
T. Scanlon, \emph{Diophantine geometry of the torsion of a Drinfeld module}.  J. Number Theory  {\bf 97}  (2002),  no. 1, 10--25.  

\bibitem{Szpiro-Ullmo-Zhang}
L. Szpiro, E. Ullmo and S. Zhang, \emph{\'{E}quir\'{e}partition des petits points}. (French) [Uniform distribution of small points]. Invent. Math. {\bf 127} (1997), no. 2, 337--347.

\bibitem{Vojta}
P. Vojta, \emph{Integral points on subvarieties of semiabelian varieties. I}. Invent. Math. {\bf 126} (1996), no. 1, 133-181.

\bibitem{Wan}
J. T.-Y. Wang, \emph{The Mordell-Weil theorems for Drinfeld modules over 
finitely generated function fields}, Manuscripta math. \textbf{106} (2001), 305--314.

\bibitem{Yafaev}
A. Yafaev, \emph{A conjecture of Yves Andr\'{e}'s}. Duke Math. J. {\bf 132} (2006), no. 3, 393--407.

\bibitem{Zhang}
S. Zhang, \emph{Equidistribution of small points on abelian varieties}. Ann. of Math. (2) {\bf 147} (1998), no. 1, 159--165.



\end{thebibliography}
\end{document}